\documentclass[11pt,oneside, final]{amsart}
\copyrightinfo{0}{Iranian Mathematical Society}
\pagespan{1}{\pageref*{LastPage}}
\usepackage{etoolbox,lastpage}
\commby{}
\usepackage{amsmath,amsthm,amscd,amsfonts,amssymb,enumerate}
\usepackage{graphicx}
\usepackage{color}
\usepackage[colorlinks]{hyperref}

\makeatletter \@addtoreset{equation}{section}

\newtheorem{theorem}{Theorem}[section]

\newtheorem{lemma}[theorem]{Lemma}
\newtheorem{corollary}[theorem]{Corollary}
\newtheorem{conjecture}[theorem]{Conjecture}

\theoremstyle{definition}

\newtheorem{remark}[theorem]{Remark}
\newtheorem{assumption}[theorem]{Assumption}
\theoremstyle{remark}
\numberwithin{equation}{section}
\newcommand{\citep}[1]{\cite{#1}}

\begin{document}


\title[Michael-Simon type inequalities in hyperbolic space]{Michael-Simon type inequalities in hyperbolic space $\mathbb{H}^{n+1}$ via Brendle-Guan-Li's flows}

\author[J. Cui]{Jingshi Cui}

\author[P. Zhao]{Peibiao Zhao$^*$}

\thanks{$^*$Corresponding author}
%
\begin{abstract}
In the present paper, we first establish and verify a new sharp hyperbolic version of the Michael-Simon inequality for mean curvatures in hyperbolic space $\mathbb{H}^{n+1}$ based on the locally constrained inverse curvature flow introduced by Brendle, Guan and Li \cite{BGL} as follows
\begin{align}\label{0.1}
    \int_{M} \lambda^{'} \sqrt{f^{2} E_{1}^{2}+|\nabla^{M} f|^{2}} -\int_{M}\left\langle \bar{\nabla}\left(f\lambda^{'}\right),\nu \right\rangle  +\int_{\partial M}  f
    \geq
    \omega_{n}^{\frac{1}{n}}\left(\int_{M}f^{\frac{n}{n-1}}\right)^{\frac{n-1}{n}}
\end{align}
provided that $M$ is $h$-convex and $f$ is a positive smooth function, where $\lambda^{'}(r)=\rm{cosh}$$r$.

In particular, when $f$ is of constant, (\ref{0.1}) coincides with the Minkowski type inequality stated by  Brendle, Hung, and Wang in \cite{BHW16}.

Further, we also establish and confirm a new sharp Michael-Simon inequality for the $k$-th mean curvatures in $\mathbb{H}^{n+1}$ by virtue of the Brendle-Guan-Li's flow  \cite{BGL} as below
\begin{align}\label{0.2}
    \int_{M} \lambda^{'} \sqrt{f^{2}E_{k}^{2} + |\nabla^{M} f|^{2} E_{k-1}^{2}} &-\int_{M}\left\langle \bar{\nabla}\left(f\lambda^{'}\right),\nu \right\rangle \cdot E_{k-1} +\int_{\partial M}  f \cdot E_{k-1} \notag \\
     & \geq
    \left(p_{k}\circ q_{1}^{-1}(W_{1}(\Omega)) \right)^{\frac{1}{n-k+1}}\left(\int_{M}f^{\frac{n-k+1}{n-k}}\cdot E_{k-1}\right)^{\frac{n-k}{n-k+1}}
\end{align}
provided that $M$ is $h$-convex and  $\Omega$ is the domain enclosed by $M$, $p_{k}(r)=\omega_{n}(\lambda^{'})^{k-1}$, $W_{1}(\Omega)  =\frac{1}{n}|M|$, $\lambda^{'}(r)=\rm{cosh}$$r$, $q_{1}(r)=W_{1}(S^{n+1}_{r})$, the area for a geodesic sphere of radius $r$, and $q_{1}^{-1}$ is the inverse function of $q_{1}$.

In particular, when $f$ is of constant and $k$ is odd, (\ref{0.2}) is exactly  the weighted Alexandrov-Fenchel inequalities proven by Hu, Li, and Wei in \cite{HLW20}.
\end{abstract}

\maketitle
%

\section{Introduction}

\noindent  The Michael-Simon inequality on the generalized
submanifolds immersed into Euclidean space was initially proposed by Michael and Simon \cite{MS1973}, and  the existence of the best constant is still an open issue. The sharp form of inequalities plays a crucial role in enhancing our understanding of the properties and applications of inequalities, such as \cite{L1994,LL2012,LL2013}. An important development is  that S. Brendle \cite{B19} in 2019 studied and  confirmed a sharp form of Michael-Simon type inequality in $\mathbb{R}^{n+1}$, where the optimal constant is the isometric constant in Euclidean space $\mathbb{R}^{n}$.
\begin{theorem}{(\cite{B19})}\label{th1.1}
    Let $M$ be a compact hypersurface in $\mathbb{R}^{n+1}$ (possibly with boundary $\partial M$), and let $f$ be a positive smooth function on $M$. Then
    \begin{align}\label{1.1}
       \int_{M} \sqrt{\left|\nabla^{M} f\right|^{2}+f^{2} H^{2}}+\int_{\partial M} f \geq n|B^{n}|^{\frac{1}{n}}\left(\int_{M} f^{\frac{n}{n-1}}\right)^{\frac{n-1}{n}},
    \end{align}
    where $H$ is the mean curvature of $M$ and $B^{n}$ is the open unit ball in $\mathbb{R}^{n}$. Moreover, if the equality holds, then $f$ is constant and $M$ is a flat disk.
\end{theorem}
In addition to the Michael-Simon inequality for mean curvatures, Sun-Yung Alice Chang and Yi Wang \cite{CW2013} introduced and derived the Michael-Simon inequality for the (non-normalized) $k$-th mean curvature $\sigma_{k}(\kappa)$ in $\mathbb{R}^{n+1}$.
\begin{theorem}{(\cite{CW2013})}\label{th1.2}
    Let $i: M^{n} \rightarrow \mathbb{R}^{n+1}$ be an isometric immersion. Let $U$ be an open subset of $M$ and $u \in C_{c}^{\infty}(U)$ be a nonnegative function. For $m =2, \cdots, n-1$, if $M$ is $(m+1)$-convex, then there exists a constant $C$ depending only on $n$, $m$ and $k$, such that for $1 \leqslant k \leqslant m$
    \begin{align}\label{1.2}
        &\left(\int_{M} \sigma_{k-1}(\kappa) u^{\frac{n-k+1}{n-k}} d \mu_{M}\right)^{\frac{n-k}{n-k+1}} \leq C \int_{M}\left(\sigma_{k}(\kappa) u+\sigma_{k-1}(\kappa)|\nabla u|+\cdots+ |\nabla^{k} u|\right) d \mu_{M}.
    \end{align}
    If $m=n$, then the inequality holds when $M$ is $n$-convex. If $m=1$, then the inequality holds when $M$ is 1-convex. ($m=1$ case is a corollary of the Michael-Simon inequality).
\end{theorem}
However, the best constant $C$ in  (\ref{1.2}) is not an absolute constant. Fortunately, Cui and Zhao  in  \cite{CZ21} developed and proved a sharp Michael-Simon inequality for the (non-normalized) $k$-th mean curvatures of hypersurfaces in $\mathbb{R}^{n+1}$.

The Michael-Simon type inequality is  a basic tool with extensive and  significant applications, especially in the arguments of the regularity of partial differential equations \cite{Ca10}, the regularity of surfaces with prescribed curvatures \cite{BG73,DHT10}, the theory of curvature flows \cite{ES92}, and so on. Therefore, it is necessary and interesting to generalize the Michael-Simon type inequality to
the case of Riemannian manifolds.

A very important work that needs to be mentioned here is that S. Brendle \cite{B21} obtained a remarkable  Michael-Simon inequality for
mean curvatures in Riemannian manifolds with nonnegative sectional curvatures  as
below.
\begin{theorem}{(\cite{B21})\label{th1.3}}
    Let $N$ be a complete noncompact manifold of dimension $n+1$ with nonnegative sectional curvatures. Let $M$ be a compact submanifold of $N$ of dimension $n$ (possibly with boundary $\partial M$), and let $f$ be a positive smooth function on $M$. Then
    \begin{align}\label{1.3}
        \int_{M} \sqrt{\left|\nabla^{M} f\right|^{2}+f^{2} |H|^{2}}+\int_{\partial M} f \geq n|B^{n}|^{\frac{1}{n}}\theta^{\frac{1}{n}}\left(\int_{M} f^{\frac{n}{n-1}}\right)^{\frac{n-1}{n}},
    \end{align}
    where $\theta$ denotes the asymptotic volume ratio of $N$ and $H$ denotes the mean curvature vector of $M$. If the equality holds, then $N$ is isometric to Euclidean space, $M$ is a flat ball, and $f$ is constant.
\end{theorem}
\begin{remark}\label{rem1.4}
    If $M$ is of closed, then the  (\ref{1.3}) in \cite{B21} is reduced naturally into the following
    \begin{align}
        \int_{M} \sqrt{\left|\nabla^{M} f\right|^{2}+f^{2} |H|^{2}} \geq n|B^{n}|^{\frac{1}{n}}\theta^{\frac{1}{n}}\left(\int_{M} f^{\frac{n}{n-1}}\right)^{\frac{n-1}{n}} \tag{1.3$^{ '}$}.
    \end{align}
\end{remark}
To this day, as far as we know,  the Michael-Simon type inequality in Riemannian manifolds with negative sectional curvatures is still open.  A counterexample in  \cite{CZ2022} shows that the Michael-Simon inequality (\ref{1.3}) does not hold in Riemannian manifolds with negative sectional curvatures.
One exciting thing is that Cui and Zhao  in \cite{CZ2022} proposed a  hyperbolic version of the Michael-Simon type inequality with respect to the $k$-th mean curvatures, which is the following conjecture for the Michael-Simon type inequality for the $k$-th mean curvatures inequality in hyperbolic space $\mathbb{H}^{n+1}$.
\begin{conjecture}\cite{CZ2022}\label{con1.5}
    Let $M$ be a compact hypersurface in $\mathbb{H}^{n+1}$ (possibly with boundary $\partial M$) and let $f$ be a positive smooth function on $M$. For $1\leq k\leq n$, there holds
    \begin{align}
        \int_{M} \lambda^{'} \sqrt{f^{2}E_{k}^{2} +|\nabla^{M} f|^{2} E_{k-1}^{2}} -\int_{M}\left\langle \bar{\nabla}\left(f\lambda^{'}\right),\nu \right\rangle \cdot |E_{k-1}| &+\int_{\partial M}  f \cdot|E_{k-1}| \notag \\
        &\geq
        C\left(\int_{M} f^{\frac{n-k+1}{n-k}}\cdot |E_{k-1}|\right)^{\frac{n-k}{n-k+1}}, \label{1.4}
    \end{align}
    and when $k=1$, we have
    \begin{align}
        \int_{M} \lambda^{'} \sqrt{f^{2} E_{1}^{2}+|\nabla^{M} f|^{2}} -\int_{M}\left\langle \bar{\nabla}\left(f\lambda^{'}\right),\nu \right\rangle  +\int_{\partial M}  f
        \geq
        C\left(\int_{M}f^{\frac{n}{n-1}}\right)^{\frac{n-1}{n}}, \label{1.5}
    \end{align}
    where $E_{k}:=E_{k}(\kappa)$ and $\nu$ are the (normalized) $k$-th mean curvature and the unit outward normal of $M$,  respectively, $\bar{\nabla}$ is the Levi-Civita connection with respect to the metric $\bar{g}$ on hyperbolic space $\mathbb{H}^{n+1}$.
\end{conjecture}

\begin{remark}\label{rem1.6}
    Conjecture \ref{con1.5} is partially resolved after adding appropriate constraints to $f$ (one can see \cite{CZ2022} for details).

    (1) By using a new locally constrained mean curvature flow, it is proved that the inequality (\ref{1.4}) holds when $M$ is starshaped.

    (2) When $M$ is starshaped and strictly $k$-convex, the inequality (\ref{1.5}) is exactly the corresponding result by virtue of Scheuer-Xia's flow \cite{SX19}.
\end{remark}
The conclusions in  \cite{CZ2022} are not particularly excellent in terms of applications, that is, it is difficult to extract or infer certain well-known outcomes or connections from the conclusions in \cite{CZ2022}. They cannot, for example, contain the scenario where $f$ is constant, which is required for the application of the Michael-Simon type inequality. Hence, the main goal of this paper is to overcome this shortcoming and thoroughly resolve Conjecture \ref{con1.5} utilizing different geometric flows.

The basic idea of curvature flow is to change the shape of hypersurfaces by continuously adjusting the geometric properties of the points on them. If the curvature flow is designed to keep one geometric quantity unchanged while another geometric quantity is monotonic along the curvature flow, then combining the asymptotic behavior of this flow one can create or prove geometric inequalities, one can refer to \cite{BHW16,CS2022,Hui87,HI2001}.

There are numerous great results concerning locally constrained flows defined mainly by the Minkowski type formula, for instance, \cite{JX19,LS20,SWX18,WX19,WX,GL2015,GL19}. Brendle, Guan and Li \cite{BGL} designed an locally constrained flow of inverse curvature type in space forms, which has an interesting feature: all quermassintegrals along the flow either have monotonicity or remain unchanged. Using this type of geometric flows, Hu, Li, and Wei \cite{HLW20} proved sharp Alexandrov-Fenchel type inequalities in $\mathbb{H}^{n+1}$ and new geometric inequalities containing the weighted integral of the $k$-th shifted mean curvature. Hu and Li \cite{HL21} also established new sharp inequalities including the weighted curvature integrals in $\mathbb{H}^{n+1}$.

In the present paper we shall use Brendle-Guan-Li's flow in \cite{BGL} to verify Conjecture \ref{con1.5}.

Let $X_{0}: \Sigma \to \mathbb{H}^{n+1}$ be a smooth embedding map such that $\Sigma_{0}=X_{0}(\Sigma)$ is a closed, starshaped and $l$-convex hypersurface. Brendle, Guan, and Li \cite{BGL} developed the locally constrained inverse curvature flow as follows:

Let  $X: \Sigma \times [0,T) \to \mathbb{H}^{n+1}$ satisfy
\begin{equation}\label{1.6}
    \begin{cases}
        \frac{\partial }{\partial t}X(x,t)=\left(\lambda^{'}\frac{E_{l-1}}{E_{l}}- u\right)\nu(x,t), \qquad l=1, \cdots, n\\
        X(\cdot,0)=X_{0}(\cdot)
     \end{cases}
\end{equation}
where $u=\left\langle \lambda\partial_{r}, \nu \right\rangle $ and $\nu$ are the support function and the unit outer normal vector of $\Sigma_{t}=X(M,t)$,  respectively. The following result has been demonstrated by Hu, Li and Wei \cite{HLW20}.

\begin{theorem}\cite{HLW20}\label{th1.7}
    Let $X_{0}: \Sigma \to \mathbb{H}^{n+1}(n\geq 2)$ be a smooth embedding of a closed, $h$-convex hypersurface $\Sigma_{0}=X_{0}(\Sigma)$ in $\mathbb{H}^{n+1}$. Then the flow (\ref{1.6}) has a unique smooth solution $\Sigma_{t}=X_{t}(\Sigma)$ for $t\in[0,+\infty)$. Moreover, $\Sigma_{t}$ is strictly $h$-convex for each $t> 0$ and it converges exponentialy to a geodesic sphere $S_{R}$ centered at the origin in $C^{\infty}$-topology as $t \to +\infty$, where the radius $R$ is uniquely determined by $W_{l}(\Omega_{0})=W_{l}(S^{n+1}_{R})$.
\end{theorem}
It is well known that if the evolution hypersurface $\Sigma_{t}$ is smooth, closed and starshaped, then the $\Sigma_{t}$  can be parameterized over $\mathbb{S}^{n}$ by the distance function $r_{t}=r(\xi,t) : \mathbb{S}^{n} \times [0, \infty) \to \mathbb{R}^{+}$ below
\begin{align*}
    \Sigma_{t}=\left\{\left(r_{t}(\xi),\xi\right)\in \mathbb{R}^{+}\times\mathbb{S}^{n}|\xi \in\mathbb{S}^{n}\right\}
\end{align*}
Moreover $\Sigma_{0}=\left\{\left(r_{0}(\xi),\xi\right)\in \mathbb{R}^{+}\times\mathbb{S}^{n}|\xi \in\mathbb{S}^{n}\right\}$, where $r_{0}=r(\xi, 0)$.

Let $M= \Sigma_{0}$ ($\Sigma_{0} \cup \partial \Sigma_{0}$ )  be a smooth starshaped and compact hypersurface (possibly with boundary), and $\forall x \in M$, we have $x=(\rho(\xi), \xi)$, which $\rho$ is the distance function of $M$
\begin{equation*}
    \rho(\xi) =
    \begin{cases}
       r_{0}(\xi),\qquad x \in \Sigma_{0}
       \\
       \bar{r}_{0}(\xi), \qquad x \in \partial \Sigma_{0}
    \end{cases}
\end{equation*}
If $f\in C^{\infty}(M)$, we have $f(x) = f (\rho (\xi), \xi)$. Also, there exists a function $\phi: \mathbb{S}^{n} \to \mathbb{R}^{+} $ such that
\begin{align*}
    f(\rho(\xi), \xi)= \phi (\xi)
\end{align*}

To prove the Michael-Simon type inequality for mean curvatures in $\mathbb{H}^{n+1}$, we use the convergence result of flow (\ref{1.6}) with $l=1$ in Theorem \ref{th1.7}, as long as $f(x)$ meets the following assumption.
\begin{assumption}\label{as1.8}
    $f(x) = \Phi_{1}\circ \rho (\xi)$ and  $\Phi_{1}$ is non-decreasing with respect to $\rho$.
\end{assumption}

\begin{theorem}\label{th1.9}
    Let $M$ be a smooth compact $h$-convex hypersurface in $\mathbb{H}^{n+1}$ (possibly with boundary $\partial M$),  the positive function $f \in C^{\infty}(M)$ satisfies Assumption \ref{as1.8}, there holds
    \begin{align}\label{1.7}
        \int_{M} \lambda^{'} \sqrt{f^{2} E_{1}^{2}+|\nabla^{M} f|^{2}} -\int_{M}\left\langle \bar{\nabla}\left(f\lambda^{'}\right),\nu \right\rangle  +\int_{\partial M}  f
        \geq
        \omega_{n}^{\frac{1}{n}}\left(\int_{M}f^{\frac{n}{n-1}}\right)^{\frac{n-1}{n}},
    \end{align}
    and when $M$ is closed, we have
    \begin{align}\label{1.8}
        \int_{M} \lambda^{'} \sqrt{f^{2} E_{1}^{2}+|\nabla^{M} f|^{2}} -\int_{M}\left\langle \bar{\nabla}\left(f\lambda^{'}\right),\nu \right\rangle
        \geq
        \omega_{n}^{\frac{1}{n}}\left(\int_{M}f^{\frac{n}{n-1}}\right)^{\frac{n-1}{n}}, \tag{1.7 ${'}$}
    \end{align}
    where $\omega_{n}$ is the area of the unit sphere $\mathbb{S}^{n}$. Equality holds in (\ref{1.8}) if and only if $M$ is a geodesic sphere.
\end{theorem}
If $f$ is of constant and $M$ is closed and $h$-convex, we can obtain the following Minkowski type inequality which was initially shown by Brendle, Hung,and Wang \cite{BHW16} for starshaped and mean convex hypersurfaces.
\begin{corollary}\label{cor1.10}
    Let $M$ be a smooth closed $h$-convex hypersurface in $\mathbb{H}^{n+1}$. Then
    \begin{align}\label{1.9}
        \int_{M}(\lambda^{'} E_{1} -u)d\mu \geq \omega_{n}^{\frac{1}{n}}|M|^{\frac{n-1}{n}},
    \end{align}
    where $u$ and $|M|$ are the support function  and the area of $M$ respectively.  Equality holds in (\ref{1.9}) if and only if $M$ is a geodesic sphere.
\end{corollary}

To prove the Michael-Simon type inequality for the $k$-th mean curvatures, $k\geq 2$, we will adopt the similar argument and  adjust the assumptions about $f(x)$.
\begin{assumption}\label{as1.11}
    (1) $f(x) = \Phi_{2}\circ \rho(\xi)$, $\Phi_{2}$ is non-decreasing with respect to $\rho$.
    \\
    (2) There exists a smooth positive function $Q(\xi,t)$ and $\frac{\partial Q}{\partial r_{t}} = \frac{\partial Q}{\partial \lambda^{'}(\xi, t)} \lambda(\xi, t)\geq 0$, such that $\widetilde{\Phi}_{2} := \Phi_{2}^{\frac{n-k+1}{n-k}}$ satisfying
    \begin{eqnarray}\label{1.10}
        \begin{cases}
            \widetilde{\Phi}_{2} - \frac{1}{n}\Delta^{\Sigma_{t}} \widetilde{\Phi}_{2} = Q (\xi,t)E_{1}(\xi,t) \qquad \text{ on } \Sigma_{t},   t\geq 0 \\
            Q(\xi , 0) = \frac{\widetilde{\Phi}_{2}}{E_{1}(\xi, 0)}
        \end{cases}
    \end{eqnarray}
\end{assumption}
\begin{remark}\label{rm1.12}
    The solution of (\ref{1.10}) exists according to the standard theory of ordinary differential equations, and $\widetilde{\Phi}_{2}$ can take constants.
\end{remark}

\begin{theorem}\label{th1.13}
    Let $M$ be a smooth, compact and  $h$-convex hypersurface in $\mathbb{H}^{n+1}$ (possibly with boundary $\partial M$), and $\Omega$ be the domain enclosed by $M$. Assume that $f$ satisfies Assumption \ref{as1.11}. Then for any $2\leq k \leq n$, there holds
    \begin{align}\label{1.11}
        \int_{M} \lambda^{'} \sqrt{f^{2}E_{k}^{2} +|\nabla^{M} f|^{2} E_{k-1}^{2}} &-\int_{M}\left\langle \bar{\nabla}\left(f\lambda^{'}\right),\nu \right\rangle \cdot E_{k-1} +\int_{\partial M}  f \cdot E_{k-1} \notag \\
        &\geq
        \left(p_{k}\circ q_{1}^{-1}\left(W_{1}(\Omega)\right) \right)^{\frac{1}{n-k+1}}\left(\int_{M}f^{\frac{n-k+1}{n-k}}\cdot E_{k-1}\right)^{\frac{n-k}{n-k+1}},
    \end{align}
    and when $M$ is closed, we have
    \begin{align}\label{1.12}
        \int_{M} \lambda^{'} \sqrt{f^{2}E_{k}^{2} +|\nabla^{M} f|^{2} E_{k-1}^{2}} &-\int_{M}\left\langle \bar{\nabla}\left(f\lambda^{'}\right),\nu \right\rangle \cdot E_{k-1} \notag \\
        &\geq
        \left(p_{k}\circ q_{1}^{-1}\left(W_{1}(\Omega)\right) \right)^{\frac{1}{n-k+1}}\left(\int_{M}f^{\frac{n-k+1}{n-k}}\cdot E_{k-1}\right)^{\frac{n-k}{n-k+1}}, \tag{1.10 $^{'}$}
    \end{align}
    where $W_{1}(\Omega)= \frac{|M|}{n}$, $p_{k}(r)=\omega_{n}(\lambda^{'})^{k-1}(r)$, $q_{1}(R)=W_{1}(S^{n+1}_{R})$ and $q_{1}^{-1}$ is the inverse function of $q_{1}$. Equality holds in (\ref{1.12}) if and only if $M$ is a geodesic sphere.
\end{theorem}

If $f=const$ in the inequality (\ref{1.12}), we can derive the weighted geometric inequality established in \cite{HLW20} between the weighted curvature integrals $W^{\lambda^{'}}_{k+1}(\Omega) =\int_{M} \lambda^{'}E_{k} d\mu $ and $W_{1}(\Omega)$ for $h$-convex domain $\Omega$.
\begin{corollary}\label{cor1.14}
    Let $\Omega$ be a $h$-convex domain with smooth boundary $M$ in $\mathbb{H}^{n+1}$. For $k=2m+1$, $m\geq 1$, there holds
    \begin{align}\label{1.13}
        W^{\lambda^{'}}_{k+1}(\Omega) \geq h_{k} \circ q_{1}^{-1}(W_{1}(\Omega)),
    \end{align}
    where $h_{k}(r)=\omega_{n}(\lambda^{'})^{k+1} \lambda^{n-k}(r)$. Equality holds in (\ref{1.13}) if and only if $\Omega$ is a geodesic ball centered at the origin.
\end{corollary}

The paper is organized as follows. Section 2 recalls some geometries of starshaped hypersurfaces in hyperbolic space $\mathbb{H}^{n+1}$, and collects some basic properties of the normalized elementary symmetric functions and derives the evolution equations. The sharp Michael-Simon type inequalities in Theorem \ref{th1.9} and Theorem \ref{th1.13} will be proved in Section 3.
\section{ Preliminaries}

\noindent In this section, we recount the known definitions and results about starshaped hypersurfaces, the normalized elementary symmetry functions, and the evolution equations along the curvature flow (\ref{2.8}).

\subsection{ Starshaped hypersurfaces in a hyperbolic
space}
\
\vglue-10pt
 \indent

The hyperbolic space $\mathbb{H}^{n+1}$ can be viewed as a warped product manifold $\mathbb{R}^{+} \times \mathbb{S}^{n}$ with the metric $\bar{g} =dr^{2}+\lambda(r)^{2}\sigma$,  where $\sigma$ is the standard metric on $\mathbb{S}^{n}\subset \mathbb{R}^{n+1}$.

Let  $\Sigma \subset \mathbb{H}^{n+1}$ be a smooth, closed starshaped hypersurface with respect to the origin, the support function $u = \left\langle \lambda\partial r, \nu \right\rangle >0$ all over on $\Sigma$, where $\nu$ is the unit outward normal vector of $\Sigma$. Furthermore, $\Sigma$ can be expressed as a radial graph over $\mathbb{S}^{n}$
\begin{align*}
    \Sigma =\{(r(\xi),\xi): \xi \in \mathbb{S}^{n}\}.
\end{align*}
Let $\xi=(\xi^{1},\cdots,\xi^{n})$ be a local coordinate of $\mathbb{S}^{n}$, $D$ be the Levi-Civita connection on $(\mathbb{S}^{n},\sigma)$, $\partial_{i}=\partial_{{\xi}^{i}}$, and $r_{i}=D_{i}r$. We introduce a new function $\varphi:\mathbb{S}^{n} \to \mathbb{R}$ by
\begin{align*}
    \varphi(\xi)=\chi(r(\xi)),
\end{align*}
where $\chi$ is a positive smooth function that satisfies $\frac{\partial }{\partial r}\chi=\frac{1}{\lambda(r)}$. Hence
\begin{align*}
    \varphi_{i}:=D_{i}\varphi=\frac{r_{i}}{\lambda(r)}.
\end{align*}
The induced metric $g_{ij}$ on $\Sigma$ and its inverse matrix $g^{ij}$ are as below
\begin{align*}
    &g_{i j}=\lambda^{2}\sigma_{i j}+r_{i} r_{j}= \lambda^{2}\left(\sigma_{i j}+ \varphi_{i} \varphi_{j}\right)  ,\\
    &g^{i j}=\frac{1}{\lambda^{2}}\left(\sigma^{ij}-\frac{r^{i} r^{j}}{\lambda^{2}+|Dr|^{2}}\right)=\frac{1}{\lambda^{2}}\left(\sigma^{ij}-\frac{\varphi^{i} \varphi^{j}}{v^{2}}\right) ,
\end{align*}
where $r^{i}=\sigma^{ik}r_{k}$, $\varphi^{i}=\sigma^{ik}\varphi_{k}$ and $v=\sqrt{1+\lambda^{-2}|Dr|^{2}}=\sqrt{1+|D\varphi|^{2}}$.
The unit outer normal $\nu$ and the support function $u$ have the following forms
\begin{align*}
    \nu= \frac{1}{v}\left(\partial_{r} - \frac{r_{i}}{\lambda^{2}}\partial_{i}\right) =\frac{1}{v}\left(\partial_{r} - \frac{\varphi_{i}}{\lambda}\partial_{i}\right),\qquad
    u=\left\langle \lambda\partial_{r}, \nu \right\rangle =\frac{\lambda}{v}.
\end{align*}
The second fundamental form $h_{ij}$, the Weingarten matrix $\mathcal{W}=(h^{i}_{j})=(g^{ik}h_{kj})$ and the mean curvature of $\Sigma$ can be expressed as  (see e.g., \cite{G11})
\begin{align*}
    &h_{ij}=\frac{\lambda^{'}}{\lambda v}g_{ij}-\frac{\lambda}{v}\varphi_{ij} \\
    &h^{i}_{j}=g^{ik}h_{kj}= \frac{\lambda^{'}}{\lambda v}\delta^{i}_{j}-\frac{1}{\lambda v}\left(\sigma^{ik}-\frac{\varphi^{i}\varphi^{k}}{v^{2}}\right)\varphi_{kj}\\
    &H=\frac{n\lambda^{'}}{\lambda v}-\frac{1}{\lambda v}\left(\sigma^{ik}-\frac{\varphi^{i}\varphi^{k}}{v^{2}}\right)\varphi_{ki}
\end{align*}

\begin{lemma}(\cite{GL2015})\label{lem2.1}
    Let $(\Sigma,g)$ be a smooth hypersurface in $\mathbb{H}^{n+1}$. We denote
    \begin{align*}
        \Gamma (r)=\int^{r}_{0} \lambda(y)dy=\lambda^{'}(r)-1.
    \end{align*}
    Then $\Gamma|_{\Sigma}$ satisfies
    \begin{align}
        \nabla_{i}\Gamma &=\nabla_{i} \lambda^{'} =\left\langle \lambda \partial r, e_{i}\right\rangle, \label{2.1} \\
        \nabla_{i}\nabla_{j}\Gamma &= \nabla_{i}\nabla_{j}\lambda^{'}=\lambda^{'}g_{ij}-uh_{ij}, \label{2.2}
    \end{align}
    where $\left\{e_{1},\cdots,e_{n}\right\}$ is a basis of the tangent space of $\Sigma$.
\end{lemma}

\subsection{ Normalized elementary symmetric
functions }
\
\vglue-10pt
 \indent

The normalized $l$-th elementary symmetric functions for $\kappa=\left(\kappa_{1}, \cdots, \kappa_{n}\right)$ are
\begin{align*}
     E_{l}(\kappa)=\binom{n}{l}^{-1} \sigma_{l}(\kappa)=\binom{n}{l}^{-1} \sum_{1 \leq i_{1}<\ldots<i_{l} \leq n} \kappa_{i_{1}} \cdots \kappa_{i_{l}},\quad l=1, \cdots, n.
\end{align*}
and usually $E_{0}(\kappa)=1$ and $E_{l}(\kappa)=0$ for $l>n$. Regarding $\kappa$ as the eigenvalues of an $n \times n$ symmetric matrix $A=[A_{ij}] $, i.e. $\kappa=\kappa(A)=(\kappa_{1},\cdots,\kappa_{n})$, we have
\begin{align*}
    E_{l}(\kappa(A))= E_{l}(A)=\frac{(n-l)!}{n!} \delta_{i_{1} \ldots i_{l}}^{j_{1} \ldots j_{l}} A_{i_{1} j_{1}} \cdots A_{i_{l} j_{l}}, \quad l=1, \ldots, n.
\end{align*}
where $\delta_{i_{1} \ldots i_{l}}^{j_{1} \ldots j_{l}}$ is a
generalized Kronecker delta. For the further information, one can see \cite{G2013} for details.

\begin{lemma}(\cite{G2013})\label{lem2.2}
    Let $\dot{E}_{l}^{i j}=\frac{\partial E_{l}}{\partial A_{i j}}$, then we have
    \begin{align}
     \sum_{i, j} \dot{E}_{l}^{i j} g_{ij} &=l E_{l-1} \label{2.3}\\
     \sum_{i, j} \dot{E}_{l}^{i j} A_{i j} &=l E_{l} \label{2.4}\\
     \sum_{i, j} \dot{E}_{l}^{i j}\left(A^{2}\right)_{i j} &=n E_{1} E_{l}-(n-l) E_{l+1} \label{2.5}
     \end{align}
    where $\left(A^{2}\right)_{i j}=\sum_{l=1}^{n} A_{i l} A_{l j}$.
\end{lemma}
Next, we show that $E_{l}$ has the divergent free structure.
\begin{lemma}(\cite{G2013})\label{lem2.3}
    Supppose $\{e_{1}, \cdots, e_{n}\}$ is a local orthonormal fram on $\Sigma$, $A = [A_{ij}]$ is a Codazzi tensor on $\Sigma$, then for each $i$,
    \begin{align*}
        \sum_{j = 1}^{n} \left(\frac{\partial E_{l}}{\partial A_{ij}}\right)_{j}(A) = 0 .
    \end{align*}
\end{lemma}
The famous Minkowski formulas and the Newton-MacLaurin inequality in $\mathbb{H}^{n+1}$ are as follows.
\begin{lemma} (\cite{GL2015})\label{lem2.4}
    Let $\Sigma$ be a smooth closed hypersurface in $\mathbb{H}^{n+1}$. Then
    \begin{align}\label{2.6}
        \int_{\Sigma} \lambda^{'}E_{l-1}(\kappa) d\mu= \int_{\Sigma} uE_{l}(\kappa) d\mu.
    \end{align}
\end{lemma}
\begin{lemma}(\cite{G2013})\label{lem2.5}
    If $\kappa \in \Gamma_{m}^{+}$, the following inequality is called the Newton-MacLaurin inequality
    \begin{align}\label{2.7}
        E_{m+1}(\kappa) E_{l-1}(\kappa) \leq E_{l}(\kappa) E_{m}(\kappa), \quad 1 \leq l \leq m.
    \end{align}
    Equality holds if and only if $\kappa_{1}=\cdots=\kappa_{n}$.
\end{lemma}

\subsection{ Evolution equations}
\
\vglue-10pt
 \indent

Along the flow of $\Sigma_{t}=X(\Sigma,t)$
\begin{align}\label{2.8}
    \begin{cases}
        \frac{\partial}{\partial t}X(x,t)=F(x,t) \nu(x,t)
       \\
       X(\cdot,0)=X_{0}
    \end{cases}
\end{align}
where $F(x,t)$ and $\nu(x,t)$ are the velocity function and  the
unit outer normal vector of $\Sigma_{t}$ respectively, then we have the following evolution equations.
\begin{lemma}(\cite{HLX14})\label{lem2.6}
    \begin{align}
        \frac{\partial}{\partial t}g_{ij}&=2F h_{ij}\label{2.9},\\
        \frac{\partial}{\partial t}d\mu_{t}&=nE_{1}F d\mu_{t}\label{2.10},\\
        \frac{\partial}{\partial t}h_{ij}&=-\nabla_{j}\nabla_{i}F+F\left((h^{2})_{ij}+g_{ij}\right)\label{2.11},\\
        \frac{\partial}{\partial t}h^{j}_{i}&=-\nabla^{j}\nabla_{i}F-F\left((h^{2})^{j}_{i}+\delta^{j}_{i}\right)\label{2.12},\\
        \frac{\partial}{\partial t}E_{l-1}&=\frac{\partial E_{l-1}}{\partial h^{j}_{i}}\frac{\partial h^{j}_{i}}{\partial t}=\dot{E}^{ij}_{l-1}\left(-\nabla_{j}\nabla_{i}F - F(h^{2})_{ij}+Fg_{ij} \right) \label{2.13},\\
        \frac{\partial}{\partial t}W_{1}(\Omega_{t})&=\int_{\Sigma_{t}} FE_{1} d\mu_{t},   \label{2.14}
    \end{align}
    where $\nabla$ denotes the Levi-Civita connection on $(M_{t},g)$.
\end{lemma}

\section{ Michael-Simon type inequalities for $k$-th mean curvatures and their applications}

\noindent In this section, we will use the smooth convergence results of Brendle-Guan-Li's flow (\ref{1.6}) to prove Theorem \ref{th1.9} and Theorem \ref{th1.13} and their applications.

\subsection{ Sharp Michael-Simon type inequality for $k$-th mean curvatures}
\
\vglue-10pt
 \indent

We will investigate and confirm the new geometric inequalities (\ref{1.7}), (\ref{1.8}), (\ref{1.11}) and (\ref{1.12}) for $h$-convex hypersurfaces in this section.

\noindent {\bf Proof of Theorem \ref{th1.9} and Theorem \ref{th1.13}}: First, we reduce the inequalities (\ref{1.7}) and (\ref{1.11}) in briefly by scaling (See e.g., proofs of Theorem 1.13 and Theorem 1.17 in \cite{CZ2022}, respectively), and obtain the following inequality
\begin{align*}
    \int_{M}E_{k-1}f^{\frac{n-k+1}{n-k}}  d\mu \geq p_{k}\circ q_{1}^{-1}\left(W_{1}(\Omega_{0})\right), \qquad 1 \leq k \leq n.
\end{align*}
Notice that  $M= \Sigma_{0}$ or $M= \Sigma_{0} \cup \partial \Sigma_{0}$, and
\begin{align}\label{3.1}
    \int_{M}E_{k-1}f^{\frac{n-k+1}{n-k}}  d\mu  \geq  \int_{\Sigma_{0}}E_{k-1}f^{\frac{n-k+1}{n-k}}  d\mu.
\end{align}
It is only necessary to prove the following inequality
\begin{align}\label{3.2}
    \int_{\Sigma_{0}}E_{k-1}f^{\frac{n-k+1}{n-k}} d\mu \geq p_{k}\circ q_{1}^{-1}\left(W_{1}(\Omega_{0})\right), \qquad 1 \leq k \leq n.
\end{align}

Furthermore, $f(x)|_{\Sigma_{0}}=\Phi_{1} \circ \rho (\xi)|_{\Sigma_{0}} = \Phi_{1} \circ r_{0}(\xi)$ when $k=1$, and $f(x)|_{\Sigma_{0}}=\Phi_{2} \circ r_{0}(\xi)$ when $k \geq 2$. For computational convenience, we introduce the parameter $s$, $s =1,2$. Denote
\begin{eqnarray*}
    \widetilde{\Phi}_{s}(r_{0}) := \Phi_{s}^{\frac{n-k+1}{n-k}}(r_{0})=
    \begin{cases}
        &\Phi_{1}^{\frac{n}{n-1}}(r_{0}), \qquad s=1, \quad i.e. k=1. \\
        &\Phi_{2}^{\frac{n-k+1}{n-k}}(r_{0}),\qquad s=2, \quad  i.e. k\geq 2.
    \end{cases} 
\end{eqnarray*}
Then, the equation (\ref{3.2}) can be expressed as
\begin{align}\label{3.3}
    \int_{\Sigma_{0}}E_{k-1} \widetilde{\Phi}_{s}(r_{0}) d\mu \geq p_{k}\circ q_{1}^{-1}\left(W_{1}(\Omega_{0})\right).
\end{align}

Secondly, we prove that the $W_{1}(\Omega_{t})$ keeps unchanged along the flow (\ref{1.6}) with $l=1$, where $\Omega_{t}$ is the domain enclosed by $M_{t}$. Using (\ref{2.6}) and (\ref{2.14}), we have
\begin{align*}
    \frac{\partial}{\partial t}W_{1}(\Omega_{t})=\int_{\Sigma_{t}}\left(\lambda^{'}-uE_{1}\right)d\mu_{t}= 0.
\end{align*}

Finally, we demonstrate the monotonicity of $\int_{\Sigma_{t}}E_{k-1} \widetilde{\Phi}_{s}(r_{0}) d\mu_{t}$ which is the crucial point for this proof.
\begin{align*}
    \frac{\partial}{\partial t} \int_{\Sigma_{t}}  E_{k-1}f^{\frac{n-k+1}{n-k}}  d\mu_{t} &= \frac{\partial}{\partial t} \int_{\Sigma_{t}}  E_{k-1}\widetilde{\Phi}_{s} (r_{0})  d\mu_{t} \\
    & =  \int_{\Sigma_{t}} \left( \frac{\partial}{\partial t} E_{k-1}\widetilde{\Phi}_{s} (r_{0}) + \frac{\partial}{\partial t}  \widetilde{\Phi}_{s} (r_{0})  E_{k-1} \right)d\mu_{t} + \int_{\Sigma_{t}}   E_{k-1}\widetilde{\Phi}_{s} (r_{0}) \frac{\partial}{\partial t} d\mu_{t} \\
    & =\int_{\Sigma_{t}}  \frac{\partial}{\partial t} E_{k-1}\widetilde{\Phi}_{s} (r_{0}) d\mu_{t} + \int_{\Sigma_{t}}   E_{k-1}\widetilde{\Phi}_{s} (r_{0}) FH d\mu_{t}.
\end{align*}
Using (\ref{2.3}), (\ref{2.5}) and (\ref{2.13}), we have 
\begin{align*}
    \frac{\partial}{\partial t} & \int_{\Sigma_{t}}  E_{k-1} \widetilde{\Phi}_{s} (r_{0}) d\mu_{t}  \\
    & = \int_{\Sigma_{t}} \left[\dot{E}^{ij}_{k-1}\left(-\nabla^{\Sigma_{t}}_{j}\nabla^{\Sigma_{t}}_{i}F - F(h^{2})_{ij}+Fg_{ij} \right)\widetilde{\Phi}_{s}  + nE_{1}E_{k-1}\widetilde{\Phi}_{s} F \right] d\mu_{t}  \\
    & = \int_{\Sigma_{t}} \left[-\dot{E}^{ij}_{k-1} \nabla^{\Sigma_{t}}_{j}\nabla^{\Sigma_{t}}_{i}\widetilde{\Phi}_{s}  F  -  \dot{E}^{ij}_{k-1}(h^{2})_{ij}\widetilde{\Phi}_{s}  F + \dot{E}^{ij}_{k-1}g_{ij} \widetilde{\Phi}_{s}  F  +nE_{1} E_{k-1}\widetilde{\Phi}_{s}  F \right]d\mu_{t} \\
    & = \int_{\Sigma_{t}} \left[-\frac{k-1}{n} g^{ij}\nabla^{\Sigma_{t}}_{j}\nabla^{\Sigma_{t}}_{i}\widetilde{\Phi}_{s} E_{k-2} F  +  (n-k+1)E_{k}  \widetilde{\Phi}_{s}  F +  (k-1)E_{k-2} \widetilde{\Phi}_{s}  F \right]d\mu_{t} \\
    &  =  \int_{\Sigma_{t}} \left[-\frac{k-1}{n} \Delta^{\Sigma_{t}}\widetilde{\Phi}_{s} + (k-1)\widetilde{\Phi}_{s} \right]E_{k-2}F + (n-k+1)\widetilde{\Phi}_{s}E_{k} F d\mu_{t}.
\end{align*}
Along the flow (\ref{1.6}) with $l=1$, $F = \left(\frac{\lambda^{'}}{E_{1}} - u\right)$, we obtain
\begin{align}\label{3.4}
    \frac{\partial }{\partial t} & \int_{\Sigma_{t}}  E_{k-1} \widetilde{\Phi}_{s} (r_{0}) d\mu_{t}  \notag \\
    =&  \int_{\Sigma_{t}} \left[-\frac{k-1}{n} \Delta^{\Sigma_{t}}\widetilde{\Phi}_{s} + (k-1)\widetilde{\Phi}_{s} \right] E_{k-2} \left(\frac{\lambda^{'}}{E_{1}} - u\right) d\mu_{t}   \\
    &+ \int_{\Sigma_{t}} (n-k+1) \widetilde{\Phi}_{s} E_{k} \left(\lambda^{'} - uE_{1}\right) d\mu_{t} . \notag
\end{align}

We then divide the discussion into two situations, $s=1$ and $s=2$, respectively.

\textbf{Case 1}: $s=1$, i.e. $k=1$. Since the first term in (\ref{3.4}) disappears, we get 
\begin{align*}
    \frac{\partial}{\partial t}\int_{\Sigma_{t}}f^{\frac{n}{n-1}}d\mu_{t} = \int_{\Sigma_{t}} n \Phi_{1}^{\frac{n}{n-1}} \left(\lambda^{'}-uE_{1}\right)d\mu_{t}.
\end{align*}
Combining (\ref{2.2}), (\ref{2.3}) and (\ref{2.4}), we obtain
\begin{align}\label{3.5}
    \dot{E}^{ij}_{l} \nabla^{\Sigma_{t}}_{i}\nabla^{\Sigma_{t}}_{j}\lambda^{'}= \dot{E}^{ij}_{l}\left(\lambda^{'}g_{ij}-uh_{ij}\right) = l \left(\lambda^{'}E_{l-1}-uE_{l}\right),\qquad l=1,\cdots,n.
\end{align}
Also, for $s=1,2$, we have
\begin{align}\label{3.6}
    \nabla_{i}^{\Sigma_{t}} \Phi_{s}(r_{0})=\left\langle \bar{\nabla}\Phi_{s}(r_{0}), e_{i}^{\Sigma_{t}}  \right\rangle = \frac{d \Phi_{s}}{ d r_{0}} \left\langle \partial_{r}, e_{i}^{\Sigma_{t}}  \right\rangle  = \frac{d \Phi_{s}}{ d r_{0}} \lambda^{-1}(r_{t})\nabla_{i}^{\Sigma_{t}}\lambda^{'}(r_{t}),
\end{align}
where $\{e_{i}^{\Sigma_{t}}\}_{1 \leq i\leq n }$ is a basis of the tangent space of $\Sigma_{t}$. Thus, using Lemma \ref{lem2.3}, (\ref{3.5}) and (\ref{3.6}), we have 
\begin{align}\label{3.7}
    \frac{\partial}{\partial t}\int_{\Sigma_{t}}f^{\frac{n}{n-1}}d\mu_{t}
    &=
    \int_{\Sigma_{t}} n \Phi_{1}^{\frac{n}{n-1}}(r_{0})\dot{E}^{ij}_{1} \nabla^{\Sigma_{t}}_{i}\nabla^{\Sigma_{t}}_{j}\lambda^{'}(r_{t}) d\mu_{t} \notag \\
    &= \int_{\Sigma_{t}} - n \dot{E}^{ij}_{1} \nabla^{\Sigma_{t}}_{i}\Phi_{1}^{\frac{n}{n-1}}(r_{0}) \nabla^{\Sigma_{t}}_{j}\lambda^{'}(r_{t}) d\mu_{t} \notag \\
    & = \int_{\Sigma_{t}} - \frac{n}{n-1} \Phi_{1}^{\frac{1}{n-1}} \lambda^{-1}(r_{t})\frac{d\Phi_{1}}{d r_{0}} |\nabla^{\Sigma_{t}}\lambda^{'}(r_{t})|^{2} d\mu_{t}  \leq 0 ,
\end{align}
where the last inequality we apply the Assumption \ref{as1.8}.

\textbf{Case 2}: $s=2$, i.e. $k \geq 2$. Using the equation (\ref{1.10}), (\ref{3.4}) can be simplified to
\begin{align}
    \frac{\partial }{\partial t} & \int_{\Sigma_{t}}  E_{k-1} \widetilde{\Phi}_{2} (r_{0}) d\mu_{t}  \notag \\
    & =  \int_{\Sigma_{t}} (k-1) Q(\xi,t) E_{k-2} \left(\lambda^{'} - uE_{1} \right) + (n-k+1) \widetilde{\Phi}_{2} E_{k} \left(\lambda^{'} - uE_{1}\right) d\mu_{t} \notag \\
    & \leq \int_{\Sigma_{t}} (k-1)Q(\xi,t)\left(\lambda^{'}E_{k-2} - uE_{k-1}\right) + (n-k+1)\widetilde{\Phi}_{2} \left(\lambda^{'}E_{k-1} - uE_{k}\right) d\mu_{t} \label{3.8}  \\
    &  = \int_{\Sigma_{t}} -\frac{k-1}{n}\frac{\partial Q}{\partial \lambda^{'}}E_{k-2}|\nabla^{\Sigma_{t}} \lambda^{'}|^{2} d\mu_{t}  -  \int_{\Sigma_{t}} \frac{n-k+1}{n} \frac{d \widetilde{\Phi}_{2} }{d r_{0}}\lambda^{-1}E_{k-1}|\nabla^{\Sigma_{t}} \lambda^{'}|^{2} d\mu_{t} \leq 0 , \label{3.9}
\end{align}
where we use (\ref{2.7}) in (\ref{3.8}), (\ref{3.5}) and Assumption \ref{as1.11} in (\ref{3.9}).

Therefore, for $s=1,2$, we obtain
\begin{align*}
   \frac{\partial}{\partial t}\int_{\Sigma_{t}}E_{k-1} \widetilde{\Phi}_{s}d\mu_{t} \leq 0 ,
\end{align*}
i.e. for any $1\leq k \leq n$, we have 
\begin{align}\label{3.10}
    \int_{\Sigma_{0}} E_{k-1}f^{\frac{n-k+1}{n-k}} d\mu &= \int_{\Sigma_{0}} E_{k-1}\widetilde{\Phi}_{s} d\mu_{t} \notag \\
    & \geq \int_{\Sigma_{t}} E_{k-1} \widetilde{\Phi}_{s}  d\mu_{t} \geq \int_{\Sigma_{\infty}} E_{k-1}\widetilde{\Phi}_{s} d\mu_{\infty}=\int_{S_{R}}E_{k-1}\widetilde{\Phi}_{s} d\mu_{\mathbb{H}^{n+1}} ,
\end{align}
where the last equality is obtained from the convergence result in Theorem \ref{th1.7} and $S_{R}=\partial S^{n+1}_{R}$, the radius $R$ determined by $W_{1}(\Omega_{0})=W_{1}(S^{n+1}_{R})$. 

In addition, $\nabla^{\Sigma_{t}} \Phi_{s} =\frac{d \Phi_{s}}{d r_{0}}\lambda^{-1} \nabla^{\Sigma_{t}}\lambda^{'}=0$ on $\Sigma_{\infty}=S_{R}$, i.e. $\Phi_{s}$ is constant on $S_{R}$. Let $\Phi_{s}\big|_{S_{R}}=\lambda^{-(n-k)}(R)$, we have
\begin{align}\label{3.11}
    \int_{S_{R}}  E_{k-1} \widetilde{\Phi}_{s}  d\mu_{\mathbb{H}^{n+1}} &= \widetilde{\Phi}_{s} (R) \int_{S_{R}} E_{k-1}(\kappa) d\mu_{\mathbb{H}^{n+1}}= \omega_{n} (\lambda^{'})^{k-1} (R) \\
    &=p_{k}\circ q^{-1}_{1}\left(W_{1}(S^{n+1}_{R})\right). \notag
\end{align}
In particular, $ \int_{S_{R}} \widetilde{\Phi}_{1}  d\mu_{\mathbb{H}^{n+1}} = \omega_{n}$.

Since $W_{1}(\Omega_{t})$ is preserved along the flow (\ref{1.6}) with $l=1$, combining (\ref{3.10}) and (\ref{3.11}), we have 
\begin{align}\label{3.12}
    \int_{\Sigma_{0}} E_{k-1}f^{\frac{n-k+1}{n-k}} d\mu \geq p_{k}\circ q^{-1}_{1}\left(W_{1}(S^{n+1}_{R})\right) = p_{k}\circ q^{-1}_{1}\left(W_{1}(\Omega_{0})\right), \qquad 1\leq k \leq n.
\end{align}
and from (\ref{3.1}), we get
\begin{align}\label{3.13}
    \int_{M} E_{k-1}  f^{\frac{n-k+1}{n-k}} d\mu_{t} \geq \int_{\Sigma_{0}} E_{k-1}  f^{\frac{n-k+1}{n-k}} d\mu  \geq p_{k}\circ q^{-1}_{1}\left(W_{1}(\Omega_{0})\right), \qquad 1\leq k \leq n.
\end{align}
Thus, the inequalities (\ref{1.7}), (\ref{1.8}), (\ref{1.11}) and (\ref{1.12}) hold.

In the following, we prove the sufficient necessary conditions for the equalities hold in the inequalities (\ref{1.8}) and (\ref{1.12}). Clearly the equalities in (\ref{1.8}) and (\ref{1.12}) hold when $M = \Sigma_{0}$ and $\Sigma_{0}$ is a geodesic sphere. We just need to show that $\Sigma_{t}$ is a geodesic sphere when the equalities hold in (\ref{1.8}) and (\ref{1.12}). 

If a smooth closed $h$-convex hypersurface $\Sigma_{t}$ attains the equalities in (\ref{1.8}) and (\ref{1.12}), we have 
\begin{align*}
    \int_{\Sigma_{t}}E_{k-1}f^{\frac{n-k+1}{n-k}} d\mu_{t} = p_{k}\circ q^{-1}_{1}\left(W_{1}(\Omega_{t})\right),
\end{align*}
and 
\begin{align}\label{3.14}
    \frac{\partial}{\partial t}\int_{\Sigma_{t}}E_{k-1}f^{\frac{n-k+1}{n-k}} d\mu_{t}=0.
\end{align}

When $k=1$, (\ref{3.14}) implies that the equality holds in (\ref{3.7}). Thus, either $\nabla^{\Sigma_{t}} \lambda^{'}(r_{t})=0$ or $\frac{d \Phi_{1}}{d r_{0}}=0$. If $\nabla^{\Sigma_{t}} \lambda^{'}(r_{t}) = \lambda(r_{t}) \nabla^{\Sigma_{t}} r_{t} = 0 $, $r_{t}$ is constant, then $\Sigma_{t}$ is a geodesic sphere. If $\frac{d \Phi_{1}}{d r_{0}}=0$, then $f=\Phi_{1}=const$ and the equality in (\ref{1.8}) can be simplified as follows
\begin{align*}
    \int_{\Sigma_{t}}(\lambda^{'}E_{1}-u)d\mu =\omega_{n}^{\frac{1}{n}}|\Sigma_{t}|^{\frac{n-1}{n}}
\end{align*}
It follows from Theorem 1.2 in \cite{BHW16} that $\Sigma_{t}$ is a geodesic sphere, for any $t>0$.

When $k\geq 2$, From Lemma \ref{lem2.5} and (\ref{3.9}) it can be concluded that $\Sigma_{t}$ is a geodesic sphere for $t>0$.

In summary, when the equalities hold in (\ref{1.8}) and (\ref{1.11}), for each $t > 0$, $\Sigma_{t}$ is a geodesic sphere. Obviosly, $\Sigma_{0}$ is also a geodesic sphere since it is smoothly approximated by a family of geodesic spheres. Moreover, from (\ref{3.11}), we have 
\begin{align*}
    \omega_{n} (\lambda^{'})^{k-1}(R) = \int_{S_{R}} E_{k-1}f^{\frac{n-k+1}{n-k}} = \left(\frac{\lambda^{'}}{\lambda}\right)^{k-1} \lambda^{n}\omega_{n}f^{\frac{n-k+1}{n-k}}.
\end{align*}
Thus, $f \big|_{S_{R}} = \lambda^{-(n-k)}(R) = const$, where $ 1 \leq k \leq n$. This concludes the proof.
\hfill${\square}$

\subsection{ Applications }
\
\vglue-10pt
 \indent

In this section, we discuss the special cases where $f$ is of constant in inequalities (\ref{1.8}) and (\ref{1.12}), respectively.

\noindent {\bf Proof of Corollary \ref{cor1.10}}: Let $f$ be of positive constant, then $\nabla^{M} f=0$ and $\bar{\nabla} f =0$. From inequality (\ref{1.8}), we have
\begin{align*}
    f\int_{M} \left(\lambda^{'} E_{1} -u\right) d\mu \geq f \omega_{n}^{\frac{1}{n}}|M|^{\frac{n-1}{n}},
\end{align*}
i.e.
\begin{align*}
    \int_{M} \left(\lambda^{'} E_{1} -u\right) d\mu \geq  \omega_{n}^{\frac{1}{n}}|M|^{\frac{n-1}{n}}.
\end{align*}
Equality holds if and only if $M$ is a geodesic sphere.
\hfill${\square}$

\noindent {\bf Proof of Corollary \ref{cor1.14}}: Let $M=\partial \Omega$ be a smooth closed hypersurface and $f=const= \lambda^{-(n-k)}(R)$ where $R$  determined by $W_{1}(\Omega)= W_{1}(S^{n+1}_{R})$, from (\ref{1.12}) and (\ref{3.12}) we have
\begin{align*}
    \int_{M} \lambda^{'}E_{k}-uE_{k-1}  d\mu & \geq \left(\omega_{n}(\lambda^{'})^{k-1}\right)^{\frac{1}{n-k+1}}\left(\int_{S_{R}} E_{k-1}d\mu_{\mathbb{H}^{n+1}}\right)^{\frac{n-k}{n-k+1}}\\
    & = \omega_{n}(\lambda^{'})^{k-1}\lambda^{n-k}(R) = \omega_{n}(\lambda^{'})^{k-1}\lambda^{n-k}\circ q^{-1}_{1}\left(W_{1}(\Omega)\right).
\end{align*}
Equivalently
\begin{align*}
    W^{\lambda^{'}}_{k+1}(\Omega)- \omega_{n}&\left(\left(\frac{W_{1}(\Omega)}{\omega_{n}}\right)^{\frac{2(n+1)}{n(k+1)}} + \left(\frac{W_{1}(\Omega)}{\omega_{n}}\right)^{\frac{2(n-k)}{n(k+1)}}\right)^{\frac{k+1}{2}} \\
    &\geq
    W^{\lambda^{'}}_{k-1}(\Omega)
    -
    \omega_{n}\left(\left(\frac{W_{1}(\Omega)}{\omega_{n}}\right)^{\frac{2(n+1)}{n(k-1)}} + \left(\frac{W_{1}(\Omega)}{\omega_{n}}\right)^{\frac{2(n-k+1)}{n(k-1)}}\right)^{\frac{k-1}{2}}.
\end{align*}
Let $k=2m+1$, $m \geq 1$, we have
\begin{align*}
    W^{\lambda^{'}}_{2m+2}(\Omega)- \omega_{n}&\left(\left(\frac{W_{1}(\Omega)}{\omega_{n}}\right)^{\frac{n+1}{n(m+1)}} + \left(\frac{W_{1}(\Omega)}{\omega_{n}}\right)^{\frac{n-2m-1}{n(m+1)}}\right)^{m+1} \\
    &\geq \cdots \geq
    W^{\lambda^{'}}_{2}(\Omega) - \omega_{n}\left(\left(\frac{W_{1}(\Omega)}{\omega_{n}}\right)^{\frac{n+1}{n}} + \left(\frac{W_{1}(\Omega)}{\omega_{n}}\right)^{\frac{n-1}{n}}\right).
    \geq 0
\end{align*}
where we use the hyperbolic Alexandrov-Fenchel-type inequality (see Theorem 1.1 in \cite{LF16})  for starshaped and strictly mean convex hypersurfaces in the last inequality. Thus
\begin{align*}
    W^{\lambda^{'}}_{k+1}(\Omega) \geq  \omega_{n}\left(\left(\frac{W_{1}(\Omega)}{\omega_{n}}\right)^{\frac{2(n+1)}{n(k+1)}} + \left(\frac{W_{1}(\Omega)}{\omega_{n}}\right)^{\frac{2(n-k)}{n(k+1)}}\right)^{\frac{k+1}{2}} = h_{k} \circ q_{1}^{-1}\left(W_{1}(\Omega)\right).
\end{align*}
Equality holds if and only if $\Omega$ is a geodesic ball centered at the origin.
\hfill${\square}$



\end{document}